\documentclass[12pt]{article}
\usepackage{amsfonts,amssymb}
\usepackage[mathscr]{eucal}
\usepackage{amsbsy}
 \makeindex \voffset = -18mm
\input xypic
\begin{document}
\baselineskip = 5mm
\newcommand \lra {\longrightarrow}
\newcommand \ra {\rightarrow}
\newcommand \de {{\vartriangle }} 
\newcommand \st {{\sf S}} 
\newcommand \sst {{\sf SS}} 
\newcommand \ab {{\sf Ab}}
\newcommand \A {{\sf A}}
\newcommand \C {{\sf C}}
\newcommand \T {{\sf T}}
\newcommand \X {{\sf X}}
\newcommand \Sm {{\sf Sm}}
\newcommand \Sch {{\sf Sch}}
\newcommand \Ni {{\sf Ni}}
\newcommand \Sv {{\sf Sv}}
\newcommand \Ps {{\sf Ps}}
\newcommand \Spc {{\sf Spc}}
\newcommand \Spt {{\sf Spt}}
\newcommand \SH {{\sf {SH}}}
\newcommand \DM {{\sf {DM}}}
\newcommand \Ch {{\sf {Ch}}}
\newcommand \uno {{S}}
\newcommand \Le {{\mathbb {L}}}
\newcommand \sg {\Sigma }
\newtheorem{theorem}{Theorem}
\newtheorem{lemma}[theorem]{Lemma}
\newtheorem{corollary}[theorem]{Corollary}
\newtheorem{example}[theorem]{Example}
\newtheorem{proposition}[theorem]{Proposition}
\newtheorem{remark}[theorem]{Remark}
\newtheorem{definition}[theorem]{Definition}
\newtheorem{conjecture}[theorem]{Conjecture}
\newcommand \pal {\! \mid \! }
\newenvironment{pf}{\par\noindent{\em Proof}.}{\hfill\framebox(6,6)
\par\medskip}
\title{{\bf Finite dimensional objects in distinguished triangles} \\ $\, $ \\}

\author{Vladimir Guletski\v \i }

\maketitle

{\scriptsize

\tableofcontents

}

\section{Introduction}
\label{intro}

Let $\C $ be a $\mathbb Q$-linear, pseudoabelian and symmetric
monoidal category with a product $\otimes $. Let $n$ be a natural
number and let $\sg _n$ be the symmetric group of permutations of
$n$ elements. For any $X\in \C $ one can define its wedge
$X^{[n)}$ and symmetric $X^{(n]}$ powers as images of the
idempotents in $End_{\C }({X^{\otimes }}^n)$ corresponding to the
"vertical" and "horizontal" irreducible representations of $\sg
_n$ over $\mathbb Q$. These powers generalize usual wedge and
symmetric powers of vector spaces over a field of characteristic
zero. Then $X$ is called to be evenly (oddly) finite dimensional
if $X^{[n)}$ (correspondingly $X^{(n]}$) is a zero object for some
natural $n$. In general, $X$ is called to be finite dimensional if
$X\cong X_+\oplus X_-$ where $X_+$ is evenly and $X_-$ is oddly
finite dimensional.

The theory of finite dimensional Chow motives was introduced by
S.-I. Kimura in \cite{Kim}, and then considered in \cite{GP1} and
\cite{GP2}. The abstract theory was developed independently by
O'Sullivan, see \cite{AK}\footnote{see also more general concept
of a Schur functor in \cite{De}}. It turns out that finite
dimensionality is closely connected with some important problems
in algebraic geometry. In particular, if $X$ is a smooth
projective complex surface without non-trivial globally
holomorphic 2-forms, then Bloch's conjecture on Albanese kernel
\footnote{see \cite{B}, \cite{J1} and \cite{Voi} for its
formulation and motivations} for $X$ holds if and only if the
motive $M(X)$ of the surface $X$ is finite dimensional, see
\cite{GP2}. For surfaces of general type \footnote{the unknown and
hard part of Bloch's conjecture} it holds iff $M(X)$ is evenly
finite dimensional, i.e. $M(X)^{[n)}=0$ for some $n$, loc. cit.

Finite dimensional objects have good properties with respect to
their tensor products and quotients\footnote{in the sense of a
pseudoabelian category}. The motive of a smooth projective curve
over a field is finite dimensional, see \cite{Kim}, Prop. 5.10,
6.9 and Th. 4.2. It follows that motives of finite quotients of
products of curves and abelian varieties are finite dimensional.
Moreover, the motives of Fermat hypersurfaces are finite
dimensional and motivic finite dimensionality is a birational
invariant for surfaces, see \cite{GP1}. Now we'd like to ask about
further properties of finite dimensional objects and apply them to
the geometry of varieties. With this in mind we consider
Voevodsky's triangulated category $\DM _{\mathbb Q}$ of motives
over a field $k$ with coefficients in $\mathbb Q$, see
\cite{Voev}, \cite{SV} and \cite{VMW} for its construction. Here
we now can use motives of Zariski open and even non-smooth
varieties, if the ground field $k$ admits a resolution of
singularities, see \cite{Voev} and \cite{SV}. For example, let $X$
be a smooth proper scheme of finite type over $k$, $Z$ a closed
subscheme in $X$, and let $U=X-Z$. Then we have the distinguished
triangle
  $$
  M(Z)\longrightarrow M(X)\longrightarrow M^c(U)
  \longrightarrow M(Z)[1]
  $$
in $\DM _{\mathbb Q}$, where $M^c$ is a motive with compact
support, see \cite{Voev}. If $X$ is a surface, then $Z$ is just a
union of curves and points on $X$. Assume that the motive of a
curve (not necessary smooth or irreducible) is finite dimensional.
Then $\dim (M(Z))<\infty $. In view of that and of birational
invariance of motivic finite dimensionality for surfaces,
\cite{GP1}, one may ask: is it true that $M(X)$ is finite
dimensional if and only if $M(U)$ is finite dimensional? Or, more
generally\footnote{the question stated by Claudio Pedrini}: is the
full subcategory generated by finite dimensional objects a thick
triangulated subcategory in $\DM _{\mathbb Q}$\footnote{as it was
pointed out to me by Bruno Kahn and Chuck Weibel, the answer is
negative in the general setting of a pseudoabelian $\mathbb
Q$-linear tensor and triangulated category}?

This paper is based on two general ideas. The first one is due to
Uwe Jannsen and consists of the expectation of a nice filtration
on wedge (correspondingly, symmetric) powers of an object $Y$,
inserted into a dist. triangle $X\to Y\to Z\to \sg X$. It should
be determined by corresponding powers of objects $X$ and $Y$,
similarly to the filtration for a short exact sequence of locally
free sheaves of modules on a manifold, see \cite{Ha}, p. 127.
Without further assumptions, however, there arise problems to show
the required compatibilities in diagrams of distinguished
triangles related to the above filtration. The second idea is then
to use a triangulated category $\T $, which is the homotopy
category of an underlying pointed model and monoidal \footnote{(i)
for short of notation, a monoidal category is symmetric and closed
monoidal; (ii) sometimes we will use the word "tensor" instead of
the word "monoidal"} category $\C $, with the monoidal structure
on $\T $ induced by the moniodal structure on $\C $ \footnote{in
other words, the localization functor $\C \to \T =Ho(\C )$ is
monoidal}. Note that $\T $ is a simplicial homotopy category, see
\cite{Ho}, so that we naturally assume that the shift functor $\sg
$ in $\T $ is a suspension $\sg X=X\wedge S^1$ by the simplicial
circle $S^1$. In so structured category $\T $ it is easy to
control powers of vertices in distinguished triangles using
cofiber sequences in the underlying category $\C $. This second
idea takes its roots in the paper \cite{May}. Our main result then
is:

\begin{theorem}
\label{main1} Let $\T $ be a pseudoabelian, $\mathbb Q$-linear,
monoidal and triangulated category, which is, at the same time,
the homotopy category of a pointed model and monoidal category $\C
$. Assume, furthermore, that the monoidal structure on $\T $ is
induced by the moniodal structure on $\C $ and that the shift
functor in $\T $ is a simplicial suspension $\sg X=X\wedge S^1$.
Then, for any distinguished triangle
  $$
  X\lra Y\lra Z\lra \sg X
  $$
in $\T $, if $X$ and $Z$ are evenly (oddly) finite dimensional, it
follows that $Y$ is also evenly (oddly) finite dimensional.
Equivalently: if $X$ is evenly (oddly) finite dimensional and $Y$
is oddly (evenly) finite dimensional, it follows that $Z$ is oddly
(evenly) finite dimensional.
\end{theorem}

The important example of the above category $\T $ is $\mathbb
Q$-localized Voevodsky's motivic stable homotopy category, which
we denote by $\SH _{\mathbb Q}$, see \cite{Voevodsky}, \cite{VW},
\cite{Mo1} and \cite{Mo1}. The underlying category $\C $ is then
the category of motivic symmetric spectra constructed by Jardine
in \cite{Jar}. The recent result due to Morel \cite{Mo3} asserts
that, if $char(k)=0$ and $-1$ is a sum of squares in $k$, then
there exists an exact and monoidal equivalence between $\mathbb
Q$-localized categories $\SH _{\mathbb Q}$ and $\DM _{\mathbb Q}$
over the ground field $k$. Applying this and Theorem \ref{main1}
to schemes of dimension one we obtain the following generalization
of Kimura's result:

\begin{theorem}
\label{main2} Let $k$ be a field of characteristic zero and let
$X$ be an integral scheme of dimension one, separated and of
finite type over $k$. Then its motive $M(X)$, considered in
Voevodsky's category $\DM _{\mathbb Q}$, is finite dimensional.
\end{theorem}

After the publishing of the first version of this paper I was
informed that the same result as in Theorem \ref{main2} has been
independently obtained by Carlo Mazza\footnote{as far as I know,
his proof involves some relevant filtration due to C. Weibel, and
Schur functors due to P. Deligne.}.

The paper is organized as follows. For the convenience of the
reader, in the second section we recall definitions and basic
results on finite dimensional objects, known results on
triangulated categories, which are homotopy categories of model
monoidal categories, contained in \cite{Ho} and \cite{May}, and
outline that the motivic stable homotopy category is an example of
such a category. In the third section we develop a homotopy
technique to deal with finite dimensionality of vertices in dist.
triangles and show the existence of the above filtration on
$Y^{[n)}$ with graded pieces $Z^{[p)}\wedge X^{[q)}$ where $p+q=n$
(and the same for symmetric powers), and then prove Theorem
\ref{main1}. In the last section we prove Theorem \ref{main2}.

{\sc Acknowledgements.} The author wish to thank to Uwe Jannsen
for many useful conversations and suggestions on the subject of
the paper, and to Jens Hornbostel for consultations on stable
homotopy categories. Also we thank Luca Barbieri-Viale, Bruno
Kahn, Fabien Morel, Ivan Panin, Claudio Pedrini and Charles Weibel
for useful discussions of the theme. The work was done while the
author enjoyed the hospitality of the University of Regensburg.

\section{Preliminary results}

\subsection{\it Basics on finite dimensional objects}
\label{basics}

Let $\C $ be a monoidal $\mathbb Q$-linear and pseudoabelian
category with a product $\otimes :\C \times \C \to \C $ and a
unite object $\uno $. For any natural number $n$ and any object
$X$ in $\C $ by $X^{(n)}$ we denote the $n$-fold product
${X^{\otimes }}^n$ in $\C $ and set $X^{(0)}=\uno $. If $f:X\to Y$
is a morphism in $\C $ then let $f^{(n)}:X^{(n)}\to X^{(n)}$ be
the $n$-fold tensor product of $f$.

Let $n$ be a natural number and let $\sg _n$ be the group of
permutations of a finite set consisting of $n$ elements. Let also
$A=\mathbb Q\sg _n$ be the group algebra (over $\mathbb Q$) of the
group $\sg _n$. A classical result asserts that the set of all
irreducible representations of $\sg _n$ over $\mathbb Q$ is in
one-to-one correspondence with the set $P_n$ of all partitions
$\lambda $ of the integer $n$, and that there exists a finite
collection $\{ e_{\lambda }\} $ of pairwise orthogonal idempotents
in $A$, such that $\sum _{\lambda \in P_n}e_{\lambda }=1_A$, and
each $e_{\lambda }$ induces the corresponding irreducible
representation of $\sg _n$ up to an isomorphism.

For any natural $n$ and for any $X\in Ob(\C )$ let $\Gamma :A\to
End_{\C }(X^{(n)})$ be a natural homomorphism sending any $\sigma
\in \sg _n$ into its "graph", i.e. the endomorphism of $M^{(n)}$
permuting factors according to $\sigma $ and the commutativity and
associativity axioms from the definition of a monoidal category
(see, for example, \cite{DeMi} or \cite{Ho}). For any $\lambda \in
P_n$ let $d_{\lambda }$ be the graph of the idempotent $e_{\lambda
}$. Since $\sum _{\lambda \in P_n}e_{\lambda }=1$ in $A$, it
follows that $\sum _{\lambda \in P_n}d_{\lambda }=1$ in $End_{\C
}(X^{(n)})$. The category $\C $ being pseudoabelian, it follows
that $X^{(n)}$ is a direct sum of the images $im(d_{\lambda })$ of
the idempotents $d_{\lambda }$.

Let $n$ be a natural number. Let $d_n^+$ be the projector
$d_{(\lambda )}$ when $\lambda $ is the partition $(1,\dots ,1)$
of $n$, and let $d_n^-$ be the projector $d_{(\lambda )}$ when
$\lambda $ is the partition $(n)$ of $n$. In other words,
  $$
  d_n^+=\frac{1}{n!}\sum _{\sigma \in \sg _n}
  sgn(\sigma )\Gamma _{\sigma }
  $$
and
  $$
  d_n^-=\frac{1}{n!}\sum _{\sigma \in \sg _n}
  \Gamma _{\sigma }
  $$
Then let
  $$
  X^{[n)}=im(d_n^+)
  $$
and
  $$
  X^{(n]}=im(d_n^-)\; .
  $$
Note that, if $X$ is a vector space over a field of characteristic
zero, then $X^{[n)}$ is a usual wedge and $X^{(n]}$ is a usual
symmetric powers of $X$, see \cite{FH}, \S B.2.

Then we say that $X$ is evenly (oddly) finite dimensional if
$X^{[n)}=0$ (correspondingly $X_-^{(n]}=0$) for some natural
number $n$. In general, $X$ is said to be finite dimensional if it
can be decomposed into a direct sum $X=X_+\oplus X_-$, such that
$X_+^{[m)}=0$ and $X_-^{(n]}=0$.

The dimension\footnote{in the sense of S.-I. Kimura} $dim(X)$ of
an evenly (oddly) finite dimensional object $X$ is, by definition,
the smallest natural number $n$, such that $X^{[n)}\neq 0$
(correspondingly, $X^{(n]}\neq 0$). For example, the motive
$\mathbb Q(d)[2d]$ is evenly one dimensional in $\DM _{\mathbb Q}$
for any $d\in \mathbb Z$. In the mixed case:
$dim(X)=dim(X_+)+dim(X_-)$. It's clear that $dim(X\oplus
Y)=dim(X)+dim(Y)$.

The following result involves just properties of pseudoabelian
monoidal categories and the theory of representations of symmetric
groups (see \cite{Kim}):

\begin{proposition}
\label{tensprod} The tensor product of two finite dimensional
objects in $\C $ is finite dimensional and a subobject\footnote{in
the sense of a pseudoabelian category} of a finite dimensional
object is also finite dimensional. Moreover, if $X$ and $Y$ are
two purely finite dimensional objects with the same parity, then
$X\otimes Y$ is evenly finite dimensional, and if $X$ and $Y$ have
different parity, then $X\otimes Y$ is oddly finite dimensional.
\end{proposition}

\begin{pf}
See \cite{Kim}, Prop. 5.10 and Cor. 5.11
\end{pf}

\begin{remark}
{\rm Geometrically it means that the motive of a fibered product
of two varieties with finite dimensional motives is finite
dimensional, and, if $X\to Y$ is a finite cover of varieties, then
$dim (M(X))<\infty $ implies $dim (M(Y))<\infty $. In particular,
any Fermat hypersurface is motivically finite dimensional and
motivic finite dimensionality is a birational invariant for
surfaces, see \cite{GP1}, Th. 2.8.}
\end{remark}

\begin{proposition}
\label{uniqdecomp} Let $X$ be a finite dimensional object in $\C $
and let $X\cong X_+\oplus X_-\cong Y_+\oplus Y_-$ be two
decompositions of $X$ into its even and odd parts. Then it follows
that $X_+\cong Y_+$ and $X_-\cong Y_-$.
\end{proposition}

\begin{pf}
\cite{Kim}, Prop. 6.3
\end{pf}

\begin{theorem}
\label{findimcurves} The motive of a smooth projective curve $X$
over a field is finite dimensional. If $M(X)=\mathbb Q\oplus
M^1(X)\oplus \mathbb Q(1)[2]$ is the decomposition\footnote{see
\cite{Mu1} or \cite{Sch}} of $M(X)$ given by a $k$-rational point
on $X$, then $M^1(X)$ is oddly finite dimensional of dimension
$2g$, where $g$ is the genus of $X$.
\end{theorem}

\begin{pf}
See \cite{She} and \cite{Kim}. The proof essentially involves the
fact that a big enough symmetric power of a smooth projective
curve is a projective bundle over its Jacobian variety by
Riemann-Roch theorem.
\end{pf}

\begin{remark}
{\rm It follows that any abelian variety and all motives, which
can be reconstructed using motives of curves via their fibered
products and finite covers are finite dimensional.}
\end{remark}

\begin{remark}
{\rm Motivic finite dimensionality controls "phantom motives": if
$M$ is a Chow-motive\footnote{see \cite{Ma}, \cite{J3}, \cite{Mu1}
and \cite{Sch} for the definition of Chow-motives} and $\dim
(M)<\infty $, then $H^*(M)=0$ implies $M=0$, see \cite{Kim},
Propositions 7.2 and 7.5, as well as Ex. 9.2.4 in \cite{AK}.}
\end{remark}

\subsection{\it Homotopy category of a pointed model monoidal category}
\label{reinforce}

The basic tool in the below proof of Theorem \ref{main1} is a
monoidal triangulated category which can be represented as the
homotopy category of a pointed model and monoidal category. Here
we recall basics on such categories following \cite{Ho}.

Let $\C $ be a pointed model and monoidal category with a monoidal
product $\wedge :\C \times \C \to \C $ and the unite object $\uno
$. The coproduct of two objects $X$ and $Y$ in $\C $ will be
denoted by $X\vee Y$. Let $f:X\to Y$ and $f':X'\to Y'$ be two maps
in $\C $. Consider the coproduct $X\wedge Y'\coprod _{X\wedge
X'}Y\wedge X'$, that is colimit of the diagram
  $$
  \diagram
  X\wedge X' \ar[rr]^-{f\wedge 1}
  \ar[dd]^-{1\wedge f'}
  & & Y\wedge X' \\ \\
  X\wedge Y'
  \enddiagram
  $$
Then so called pushout smash product of $f$ and $f'$ is a unique
map
  $$
  f\square f':X\wedge Y'\coprod _{X\wedge X'}Y\wedge X'\lra
  Y\wedge Y'
  $$
determined by the above colimit. The connection between model and
monoidal structures can be then expressed by the following two
axioms, see \cite{Ho}, 4.2:

\begin{itemize}

\item
If $f$ and $f'$ are cofibrations then $f\square f'$ is also a
cofibration. If, in addition, one of two maps $f$ and $f'$ is a
weak equivalence, then so is $f\square f'$.

\item
If $q:Q\uno \to \uno $ is a cofibrant replacement for the unite
object $\uno $, then the maps $q\wedge 1:Q\uno \wedge X\to \uno
\wedge X$ and $1\wedge q:X\wedge Q\uno \to X\wedge \uno $ are weak
equivalences for all cofibrant $X$.

\end{itemize}

Now we need some series of sophisticated definitions. Let $\sf M$
be a category and let $\sf R$ be a monoidal category with a
product $\otimes $ and the unite object $\uno $. Then $\sf M$ is
called to be a (right) $\sf R$-module if we have a functor
$\otimes :\sf M\times \sf R\to \sf M$ and two natural coherent
isomorphisms $(M\otimes K)\otimes L\cong M\otimes (K\otimes L)$
and $M\otimes \uno \cong M$ for any $M\in Ob(\sf M)$, $K,L\in
Ob(\sf R)$. An $\sf R$-module is closed if there is the
corresponding adjunction, see \cite{McL}. If $\sf A$ and $\sf B$
are two monoidal categories and $F:\sf A\to \sf B$ is a monoidal
functor, then we say that $\sf B$ is an algebra over $\sf A$, see
\cite{Ho}, p. 104.

Let $\C $ be a monoidal model category. A model category $\sf D$
is called to be a $\C $-model category if (i) it is a right $\C
$-module, (ii) the action $\sf D\times \C \to \sf D$ is a Quillen
bifunctor\footnote{see \cite{Ho}, p. 107} and (iii) for an
cofibrant $X\in Ob(\sf D)$ and for cofibrant replacement $q:Q\uno
\to \uno $ the map $1\otimes q:X\otimes Q\uno \to X\otimes \uno $
is a weak equivalence, see \cite{Ho}, p. 114.

Let now $\de $ be the category of ordered finite sets and order
preserving maps between them. Let $\st $ be the category of sets
and let $\sst =\de ^{op}\st $ be the category of simplicial sets.
We always consider $\sst $ with a standard model structure on it.
Correspondingly, one has a pointed model category of simplicial
sets $\sst _*$. Then an $\sst $-model category is called a
simplicial model category\footnote{see an equivalent definition in
Ch. II, \S 2 and 3 of \cite{Jar2}}. A pointed simplicial model
category can be considered as a $\sst _*$-model category, see
Prop. 4.2.19 in \cite{Ho}.

The homotopy category $Ho(\C )$ of any simplicial model category
$\C $ (pointed simplicial model category $\C $) is a closed $\sst
$-module ($\sst _*$-module). Moreover, this fact can be geralized
to an arbitrary situation: the homotopy category $Ho(\C )$ of a
model category $\C $ (pointed model category $\C $) is a closed
$\sst $-module ($\sst _*$-module), see \cite{Ho}, Ch.5. If the
category $\C $ is monoidal (pointed monoidal), then the monoidal
structure on $\C $ induces a monoidal structure on $Ho(\C )$, see
\cite{Ho}, Th.4.3.2, and, moreover, $Ho(\C )$ is an algebra over
$Ho(\sst )$ (over $Ho(\sst _*)$), loc. cit., Ch.5.

Now we'd like to recall (following \cite{Ho}) a triangulated
structure on $Ho(\C )$, compatible with the monoidal structure, so
that $Ho(\C )$ is a monoidal (=tensor) triangulated category in
the sense of Appendix 8A in \cite{VMW}. Let $\C $ be a pointed
model and monoidal category. Since it is pointed, it's better to
denote its product by $\wedge $. The category $Ho(\C )$ is an
algebra over $Ho(\sst _*)$. Then one can consider smash products
$X\wedge K$ of any object $X\in Ob(Ho(\C ))=Ob(\C )$ with any
pointed simplicial set $K$, in particular, with the simplicial
interval $I$ and simplicial circle $S^1$. Therefore, for any $X$
we have a natural notion of its cone $CX=X\wedge I$ and suspension
$\sg X=X\wedge S^1$, as well as mapping cone $Cf$, say, for any
cofibration $f:X\to Y$ between cofibrant objects defined as a
colimit of the diagram
  $$
  \diagram
  X\ar[rr]^-{}\ar[dd]^{} & & CX \\ \\
  Y & &
  \enddiagram
  $$
Assume now that $\sg $ is a Quillen equivalence with an adjoint
loop functor $\Omega $. Then $Ho(\C )$ is naturally a triangulated
category with an endofunctor $\sg $, see \cite{Ho}, 6.5, 6.6, 7.1,
and the triangulated structure is compatible with the monoidal
structure in the usual sense, see \cite{VMW}, A8. To be more
precise, $Ho(\C )$ is a triangulated category in Hovey's sense,
i.e. it is a pre-triangulated category (see \cite{Ho}, 6.5) and
the suspention functor $\sg $ is an autoequivalence on $Ho(\C )$.
It can be shown that any triangulated category in the Hovey's
sense is a classically triangulated category, see \cite{Ho}, Prop.
7.1.6.

The triangulated category $\T =Ho(\C )$ has at least two important
advantages: strong compatibility of monoidal and triangulated
structures in $\T $, see \cite{May}, and the possibility to
describe distinguished triangles in $\T $ in terms of cofiber
sequences in $\C $, see \cite{Ho} and \cite{May}, pp. 18-19. If
$f:X\to Y$ is a map in the category $\T $, then, using cofibre
replacement in $\C $, one can assume that $f$ is a cofibration
between cofibrant objects $X$ and $Y$. Let $Z=Y/X$ be a quotient
object, that is colimit of the diagram
  $$
  \diagram
  X\ar[rr]^-{}\ar[dd]^{f} & & Y \\ \\
  * & &
  \enddiagram
  $$
where $*$ is a initial-terminal object in the pointed category $\C
$. Then $\sg X$, being a cogroup object, coacts on $Z$, see
\cite{Ho}, Th. 6.2.1. In particular, one can define a standad
baundary map $\partial :Z\to \sg X$ as a composition of the
coaction $Z\to Z\coprod \sg X$ with the evident map $Z\coprod \sg
X\to \sg X$. Then we have a cofiber distinguished triangle
  $$
  X\stackrel{f}{\lra }Y\lra Z\stackrel{\partial }{\lra }\sg X
  $$
where $Z$ is a quotient object $Y/X$. And, in fact, any
distinguished triangle in $\T =Ho(\C )$ is isomorphic to a cofiber
dist. triangle of the above type, see \cite{Ho}, 6.2-7.1, as well
as \cite{May}, Section 5.

\begin{lemma}
\label{natur} Let $f:X\to Y$ and $f':X'\to Y'$ be two cofibrations
of cofibrant objects in $\C $ with cofibers $Z$ and $Z'$
correspondingly. Let $a:X\to X'$ and $b:Y\to Y'$ be maps, such
that $bf=f'a$, and let $c:Z\to Z'$ be the induced map on cofivers.
Then $c$ is equivariant in $\T $ with respect to the cogroup
homomorphism $\sg a$, so that we have the corresponding map of
distinguished triangles in $\T $:
  $$
  \diagram
  X \ar[rr]^-{f} \ar[dd]^-{a} & & Y \ar[rr]^-{} \ar[dd]^-{b}
  & & Z \ar[rr]^-{\partial } \ar[dd]^-{c} & & \sg X \ar[dd]^-{\sg a} \\ \\
  X' \ar[rr]^-{f'} & & Y' \ar[rr]^-{} & & Z' \ar[rr]^-{\partial } & & \sg X'
  \enddiagram
  $$
\end{lemma}

\begin{pf}
See \cite{Ho}, Prop. 6.2.5.
\end{pf}

\subsection{\it Motivic stable homotopy category}

Let's consider now an important particular case of the above
abstract situation. Let $k$ be a field and let $\Sm $ be the
category of all smooth schemes, separated and of finite type over
$k$\footnote{further called simply smooth schemes}. Let $\Ps (\Sm
)$ be the category of presheaves of sets on $\Sm $. Let, further,
$\Spc =\de ^{op}\Ps (\Sm )$ be the category of simplicial
presheaves\footnote{the homotopy categories of simplicial sheaves
and presheaves are canonically isomorphic via the forgetfull
functor, see \cite{Jar}, Th. 1.2 (2), p. 453, and we prefer to
follow notation of \cite{Jar}} of sets on $\Sm $ and let also
$\Spc _*$ be the corresponding pointed category with an evident
notion of a terminal-initial object $*$. Simplicial presheaves on
$\Sm $ play the role of spaces in algebraic topology, whence the
notation $\Spc $ and $\Spc _*$, see \cite{Voevodsky}. The model
structure on $\Spc $ (and, therefore, on $\Spc _*$) depends on the
Nisnevich topology on the category $\Sm $, see \cite{Jar3},
\cite{Jar}, \cite{Voevodsky}, \cite{Mo1} and \cite{Mo2}. The
composition of the Yoneda embedding with the functor from
presheaves into simplicial presheaves allows to identify a smooth
scheme with the corresponding simplicial presheaf, see
\cite{Voevodsky}, \cite{Mo1} and \cite{Mo2}. Then, since $\Spc $
is cocomplete, one can consider colimits of spaces, for example,
quotiens, contractions, glueings, etc. in $\Spc $. In particular,
let
  $$
  T=\mathbb A^1/(\mathbb A^1-0)
  $$
be the quotient of $\mathbb A^1$ by $\mathbb A^1-0$, i.e. the
colimit of the diagram
  $$
  \diagram
  \mathbb A^1-0\ar[rr]^-{}\ar[dd]^{} & & \mathbb A^1 \\ \\
  * & &
  \enddiagram
  $$
In the homotopy category $Ho(\Spc _*)$ one has
  $$
  T\cong \mathbb P^1\cong S^1\wedge (\mathbb A^1-0)
  $$
where $S^1$ is a simplicial circle represented by the
corresponding constant (and pointed) simplicial presheaf coming
from simplicial sets.

Now a $T$-spectrum $X$ (or a motivic spectrum) is a sequience of
objects $X^n\in \Spc _*$ and bonding maps $T\wedge X^n\to X^{n+1}$
for each $n$. A map of spectra $f:X\to Y$ cosists of maps
$f^n:X^n\to Y^n$ commuting with bonding maps. A motivic symmetric
spectrum $X$ is a motivic spectrum $X$ with an extra (left) action
of the symmetric group $\sg _n$ on each $X^n$ and with $\sg
_m\times \sg _p$-equivariant compositions of bonding maps
$T^m\wedge X^n\to X^{m+n}$, where $T^m=T^{(m)}$. A map of motivic
symmetric spectra should be, of course, also equivariant for the
symmetric group action. All of these can be found in \cite{Jar},
Section 4, and should be compared with the theory in \cite{HSS}.

Let $\Spc ^{\sg }_T$ be the category of motivic symmetric spectra.
In \cite{Jar}, Th. 4.15, Jardine constructed a model structure on
$\Spc ^{\sg }_T$ and Th. 4.30, loc. cit., asserts that the
resulting homotopy category $Ho(\Spc ^{\sg }_T)$ is the desired
motivic stable homotopy category $\SH $ on the Nisnevich site, see
\cite{Voevodsky}, \cite{Mo1}, \cite{Mo2}.

Since the category $\Spc ^{\sg }_T$ is a pointed monoidal model
category, it follows a monoidal structure on $\SH $, such that the
corresponding localization functor $\Spc ^{\sg }_T\to \SH $ is
monoidal, see Th. 4.3.2 in \cite{Ho}. Moreover, applying the above
general method, breafly described in Section \ref{reinforce} (see
\cite{Ho}, Sections 6.1 - 7.1 for detailes), one can construct the
structure of a trianglulated category on $\SH =Ho(\Spc ^{\sg }_T)$
with the shift functor being induced by the suspension $\sg X=
X\wedge S^1$. So, the category $\SH $ is a typical example of a
category in the assumptions in Theorem \ref{main1}.

It's important to connect now the motivic stable homotopy category
$\SH $ with the triangulated category $\DM $ of motives over the
ground field $k$, built in \cite{Voev}. One can construct another
monoidal and triangulated category $\tilde {\DM }$, see
\cite{Mo2}, Sections 4.3 - 5.2, of, so called, $\mathbb
P^1$-motivic unbounded complexes \footnote{or $\mathbb
P^1$-stabilized category of unbounded complexes} of abelian
sheaves on the Nisnevich site. Taking free presheaves of abelian
groups, generated by presheaves of sets, and associating
normalized chain complexes to free generated simplicial
presheaves, one can construct a monoidal triangulated functor
  $$
  C_*:\SH \lra \tilde {\DM }\; ,
  $$
which is a monoidal triangulated equivalence after the
localization by $\mathbb Q$, see \cite{Mo2}, p. 32. Moreover,
there is a canonical triangulated and monoidal functor $\Phi
:\tilde {\DM }\to \DM $, which is not an equivalence in general.
The recent result due to Fabien Morel asserts that, if $char(k)=0$
and $-1$ is a sum of squares in $k$, then $\Phi $ induces an
equivalence of categories after the localization by $\mathbb Q$:

\begin{theorem}
\label{MorelTheorem} Let $k$ be a field, such that $char(k)=0$ and
$-1$ is a sum of squares in $k$. Then
  $$
  \Phi C_*:\SH _{\mathbb Q}\lra \DM _{\mathbb Q}
  $$
is a monoidal and triangulated equivalence of categories.
\end{theorem}

\begin{pf}
\cite{Mo2}, \cite{Mo3} + \cite{Mo4}
\end{pf}

\begin{remark}
\label{adddm} {\rm If $F:\T \to \T '$ is an additive and monoidal
equivalence of $\mathbb Q$-linear monoidal categories, then $X\in
Ob(\T )$ is (evenly, oddly) finite dimensional in $\T $ iff $FX\in
Ob(\T ')$ is (evenly, oddly) finite dimensional in $\T '$. Thus we
see now that Theorem \ref{main1} implies the same additivity for
finite dimensional objects in distinguished triangles in the
category $\DM _{\mathbb Q}$, constructed over an arbitrary field
$k$ of characteristic zero\footnote{the condition $-1$ to be a sum
of squares is not important since we consider motives with
coefficients in $\mathbb Q$}.}
\end{remark}

\begin{remark}
\label{ImpRem} {\rm Let $X$ be a smooth projective complex surface
with $p_g=0$. Let $M(X)$ be its motive in $\DM _{\mathbb Q}$ and
let $E^{\sg }(X)$ be its symmetric spectrum in $\SH _{\mathbb Q}$.
It is known, see \cite{GP2}, Theorem 27, that Bloch's conjecture
holds for $X$ if and only if $M(X)$ is finite dimensional in $\DM
_{\mathbb Q}$. Applying Theorem \ref{MorelTheorem}, we have that
Bloch's conjecture holds for $X$ if and only if $E^{\sg }(X)$ is
finite dimensional in $\SH _{\mathbb Q}$. This gives two
interesting sides of the matter: a motivic stable homotopy view on
the Bloch conjecture and, at the same time, a new application of
the motivic stable homotopy category to a very concrete
geometrical problem.}
\end{remark}

\section{Finite dimensional objects in dist. triangles}

\subsection{\it Cofiber sequences and combinatorics of powers}
\label{cofibreseq}

Let $\T $ be a $\mathbb Q$-linear, monoidal and triangulated
category reinforced by a pointed closed monoidal and model
category $\C $ in the sense of Section \ref{reinforce}, so that
$\T =Ho(\C )$. The monoidal product in $\C $ we denote by $\wedge
$ and the coproduct by $\vee $. The monoidal product in $\T $ will
be denoted by $\otimes $ and the direct sum by the symbol $\oplus
$. The canonical (localization) functor $\C \to \T $ is monoidal,
i.e. it carries an object $X\wedge Y$ in $\C $ into the object
$X\otimes Y$ in $\T $. The endofunctor $X\mapsto \sg X=X\wedge
S^1$ is a suspension by a simplicial circle $S^1$. Let's also
recall that "monoidal" always means "closed and symmetric
monoidal". In particular, for any fibrant $X\in Ob(\C )$ both
functors $-\wedge X$ and $X\wedge -$ preserve colimits in $\C $.

Let
 $$
 X\stackrel{f}{\lra }Y\stackrel{g}{\lra }Z\stackrel{h}{\lra }\sg X
 $$
be a distinguished triangle in $\T $. Without loss of generality,
applying cofibrant replacement, we may assume that both $X$ and
$Y$ are cofibrant and the above dist. triangle is a cofibration
triangle. In particular, $Z=Y/X$.

For any natural number $m$ a set of the cardinal $m$ will be
called $m$-set and a $0$-set is, by definition, an empty set. Fix
a natural number $m$. For any integer $0\leq i\leq m$ let
$\mathcal S_i$ be a set of $i$-subsets in the set $\{ 1,\dots ,m\}
$. For any $S\in \mathcal S_i$ let
  $$
  (Y,X)_S
  $$
be a smash-product $A_1\wedge \dots \wedge A_m$ in $\C $, such
that, for each $j\in \{ 1,\dots ,m\} $, the object $A_j$ coincides
with $X$ if $j\in S$, and with $Y$ otherwise. Surely, there are
$C^i_m$ objects of the type $(Y,X)_S$. For any $S\in \mathcal S_i$
let
  $$
  t_S:X^{(m)}\to (Y,X)_S
  $$
be an evident morphism induced by the morphisms $1_X$ and $f$. Let
  $$
  (Y,X,m-i,i)
  $$
be a colimit in $\C $ of all morphisms $t_S$. In particular,
$(Y,X,m,0)=Y^{(m)}$ and $(X,Y,0,m)=X^{(m)}$. Let also
  $$
  t_{m,i}:X^{(m)}\to (Y,X,m-i,i)
  $$
be the corresponding canonical map.

Since smashings with cofibrants preserve colimits, for any $i\in
\{ 1,\dots ,m\} $ there exists a canonical morphism
  $$
  w_{m,i}:(Y,X,m-i,i)\lra (Y,X,m-i+1,i-1)
  $$
in $\C $ induced by the morphism $f:X\to Y$ applied on different
factors. In particular, $w_{2,1}$ is the pushout square $f\square
f$. Let
 $$
 v_{m,i}=w_{m,1}\circ \dots \circ w_{m,i}\; ,
 $$
so that $v_{m,i}$ is a map
  $$
  v_{m,i}:(Y,X,m-i,i)\lra Y^{(m)}\; .
  $$
Then $v_{m,1}$ can be considered as a multi-pushout product of $f$
by $f$. Clearly, $v_{m,m}$ coincides with the morphism
$f^{(m)}:X^{(m)}\to Y^{(m)}$. More generally, we have two
commutative diagrams
  $$
  \diagram
  (Y,X,m-i,i) \ar[rr]^-{w_{m,i}}
  \ar[ddrr]^{v_{m,i}} & & (Y,X,m-i+1,i-1)
  \ar[dd]^-{v_{m,i-1}} \\ \\
  X^{(m)} \ar[rr]^-{f^{(m)}}
  \ar[uu]^-{t_{m,i}}
  & & Y^{(m)}
  \enddiagram
  $$
and
  $$
  \diagram
  (Y,X,m-i,i) \ar[rr]^-{w_{m,i}}
  & & (Y,X,m-i+1,i-1)
  \ar[dd]^-{v_{m,i-1}} \\ \\
  X^{(m)} \ar[rr]^-{f^{(m)}}
  \ar[uu]^-{t_{m,i}} \ar[uurr]^{t_{m,i-1}}
  & & Y^{(m)}
  \enddiagram
  $$

For any $S\in \mathcal S_i$ let $(Z,X)_S$ be a smash product
$B_1\wedge \dots \wedge B_m$ where, for each $j\in \{ 1,\dots ,m\}
$, $B_j=X$ if $j\in S$ and $B_j=Z$ otherwise. Let also $[Z,X]_S$
be $(Z,X)_S$, but considered in $\T $, that is to say $[Z,X]_S$ is
a tensor product $B_1\otimes \dots \otimes B_m$ where, for each
$j\in \{ 1,\dots ,m\} $, $B_j=X$ if $j\in S$ and $B_j=Z$
otherwise. There are $C^i_m$ objects of type $[X,Z]_S$. Let
$(Z,X,m-i,i)$ be a coproduct $\vee _{S\in \mathcal S_i}(Z,X)_S$ in
the category $\C $. And let $[Z,X,m-i,i]$ be a direct sum $\oplus
_{S\in \mathcal S_i}[Z,X]_S$ in the triangulated category $\T $,
that is the object $(Z,X,m-i,i)$ viewing as an object in $\T $.

\begin{proposition}
\label{cofibdisttr} The morphism $w_{m,i}$ is a cofibration for
any $i$. Moreover, the corresponding quotient object
$(X,Y,m-i+1,i-1)/(X,Y,m-i,i)$ is canonically isomorphic in $\T $
to the object $[Z,X,m-i+1,i-1]$, so that we have a cofibration
distinguished triangle

  {\small
  $$
  (Y,X,m-i,i)\stackrel{w_{m,i}}{\ra }(Y,X,m-i+1,i-1)\ra
  [Z,X,m-i+1,i-1]\ra \sg (Y,X,m-i,i)
  $$
  }

\noindent in the category $\T $.
\end{proposition}

\begin{pf}
Since $\C $ is a closed monoidal model category, it follows that,
for any cofibrant object $D$ in $\C $, both functors $-\wedge D:\C
\to \C $ and $D\wedge -:\C \to \C $ are Quillen functors, see
\cite{Ho}, 4.2. In particular, they preserve cofibrations. By our
assumtion, $X$ and $Y$ are cofibrant, so that $f^{(m)}:X^{(m)}\to
X^{(m)}$ is a cofibration. More generally: if $T$ is a $j$-subset
in an $i$-set $S\in \mathcal S_i$, then the corresponding map
$(Y,X,m,j)\to (Y,X,m,i)$ is a cofibration. It follows that any map
$t_{m,i}$ is also a cofibration. Since $f^{(m)}$ is a cofibration
and $f^{(m)}=v_{(m,i)}\circ t_{(m,i)}$, we have that $v_{(m,i)}$
is a cofibration by 2-of-3 lemma. By the pushout product axiom for
monoidal model category, the map $w_{m,1}=v_{m,1}=f\square \dots
\square f$ is a cofibration. This is a basis for an induction.
Assume that all $w_{m,1},\dots ,w_{m,i-1}$ are cofibrations. Then
their composition $v_{m,i-1}$ is a cofibration. The map
$v_{(m,i)}$ being a cofibration it follows that $w_{(m,i)}$ is a
cofibration, again by 2-of-3 lemma.

To show the second assertion of the lemma let's recall that $\C $
is a pointed category, so that, for any two objects $A$ and $B$ in
$\C $, their coproduct $A\vee B$ is a colimit of the diagram
  $$
  \diagram
  * \ar[rr]^-{} \ar[dd]^-{} & & B \\ \\
  A
  \enddiagram
  $$
If $A\to B$ is a cofibration, then the quotient object $A/B$ is a
colimit of the diagram
  $$
  \diagram
  A \ar[rr]^-{} \ar[dd]^-{} & & B \\ \\
  *
  \enddiagram
  $$
and it may be thought as a contraction of the subobject $A$ in $B$
into the fixed point $*\to B$ in $B$. Then it is easy to see that
the quotient object $(X,Y,m-i+1,i-1)/(X,Y,m-i,i)$ is exactly the
coproduct $(Z,X,m-i,i)$ in $\C $ (for the proof one has just to
draw colimit diagrams carefully). Therefore, as an object it $\T
$, it is equal to the direct sum $[Z,X,m-i,i]$, so that we obtain
the desired cofibration exact triangle from the statment of the
lemma.
\end{pf}

For example, let $m=2$. Then
  $$
  (Y,X,1,1)=(Y\wedge X)\coprod _{X\wedge X}(X\wedge Y)
  $$
and we have the following commutative diagram
  $$
  \diagram
  (Z\otimes X)\oplus (X\otimes Z) & & & & \\ \\
  (Y\wedge X)\coprod _{X\wedge X}(X\wedge Y)
  \ar[uu]^-{} \ar[rr]^-{w_{2,1}} & &
  Y\otimes Y \ar[rr]^-{g\otimes g} & & Z\otimes Z \\ \\
  X\otimes X \ar[uu]^-{w_{2,2}} \ar[rruu]^-{f\otimes f}
  \enddiagram
  $$
where the vertical row and the horizontal column both are
cofibration dist. triangles from Prop. \ref{cofibdisttr}.

Let $m=3$. Then we have three cofiber distinguished triangles
coming from Proposition \ref{cofibdisttr}. The dist. triangle "of
the first iteration" is
  $$
  \diagram
  Z\otimes Z\otimes Z \\ \\
  Y\otimes Y\otimes Y \ar[uu]^-{g^{(3)}} \\ \\
  (X\wedge Y\wedge Y)\coprod _{X\wedge X\wedge X}
  (Y\wedge Y\wedge X)
  \coprod _{X\wedge X\wedge X}(Y\wedge X\wedge Y)
  \ar[uu]^-{w_{3,1}}
  \enddiagram
  $$
Here the map $w_{3,1}$ is a fiber (in the sense of triangulated
categories) of the map $g^{(3)}$. Or, in a homotopical language:
$Z^{(3)}$ is a quotient of $Y^{(3)}$ by a subobject $(Y,X,2,1)$.
The dist. triangle of the second iteration looks like
  $$
  \diagram
  (X\otimes Z\otimes Z)\oplus (Z\otimes Z\otimes X)
  \oplus (Z\otimes X\otimes Z) \\ \\
  (X\wedge Y\wedge Y)\coprod _{X\wedge X\wedge X}(Y\wedge Y\wedge X)
  \coprod _{X\wedge X\wedge X}(Y\wedge X\wedge Y)
  \ar[uu]^-{} \\ \\
  (X\wedge Y\wedge X)\coprod _{X\wedge X\wedge X}(Y\wedge X\wedge X)
  \coprod _{X\wedge X\wedge X}(X\wedge X\wedge Y)
  \ar[uu]^-{w_{3,2}}
  \enddiagram
  $$
and the third one is
  $$
  \diagram
  (X\otimes Z\otimes X)\oplus (Z\otimes X\otimes X)
  \oplus (X\otimes X\otimes Z)
  \\ \\
  (X\wedge Y\wedge X)\coprod _{X\wedge X\wedge X}(Y\wedge X\wedge X)
  \coprod _{X\wedge X\wedge X}(X\wedge X\wedge Y)
  \ar[uu]^-{} \\ \\
  X\otimes X\otimes X
  \ar[uu]^-{}
  \enddiagram
  $$

\subsection{$2^m$-diagram}
\label{2^m}

Here we show how to build commutative diagrams illustrating the
above iterations. Let $K_m$ be $m$-dimensional cube in $\mathbb
R^m$ and let $V$ be a set of vertices of $K_m$. Fix a vertix
$v_0\in V$ a "place" the object $Y^{(m)}$ on it. Let $V_1$ be the
set of all verices of $K_m$ adjoint with $v_0$. Place objects
$(Y,X)_S$ in vertices from $V_1$ when $S$ runs $\mathcal S_1$ and
orient edges connecting vertices from $V_1$ with $v_0$ so that
they start in $V_1$ and has $v_0$ as a target. Such oriented edges
from $V_1$ to $v_0$ can be considered as cifibrations from the
above objects $(Y,X)_S$, $S\in \mathcal S_1$, into $Y^{(m)}$
induced by the cofibration $f:X\to Y$. Let $V_0=\{ v_0\} $ and let
$V_2$ be all vertices from $V\backslash (V_1\cup V_0)$ adjoint
with vertices from $V_1$. Place objects $(Y,X)_S$, $S\in \mathcal
S_2$ in vertices from $V_2$ and orient edges connecting vertices
from $V_2$ with vertices from $V_1$ so that they start in $V_2$
and has vertices from $V_1$ as their targets. Again, these
oriented edges are the cofibrations from the objects $(Y,X)_S$,
$S\in \mathcal S_2$, into the objects $(Y,X)_S$, $S\in \mathcal
S_1$, induced by the cofibration $f$. And so on. On the final step
we obtain that in the vertix $v_{2^m-1}$ opposite to $v_0$ we
place the object $X^{(m)}$ and consider the edges from $v_{2^m-1}$
to the vertices from $V_{m-1}$ as the cofibrations from $X^{(m)}$
into the objects $(Y,X)_S$ where $S$ runs $\mathcal S_{m-1}$.

As a result we obtain the commutative $2^m$-diagram associated
with the cofibration $f:X\to Y$ and the natural number $m$. We
will denote it by the same symbol $K_m$. For any $i\in \{
0,1,\dots ,m\} $ the objects $(Y,X)_S$, $S\in \mathcal S_i$, can
be identified with vertices from $V_i$. We will say that the
objects from $V_i$ are objects of the $i$-th iteration. In this
terminology, any object $(Y,X,m-i,i)$ is a colimit of all
compositions of arrows in $K_m$ starting in $X^{(m)}$ and with a
target in the objects of the $i$-th iteration.

For example, $2^2$-diagram is
  $$
  \diagram
  X\wedge X \ar[dd]^-{f\wedge 1}
  \ar[rr]^-{1\wedge f} & &
  X\wedge Y \ar[dd]^-{f\wedge 1} \\ \\
  Y\wedge X \ar[rr]^-{1\wedge f} & & Y\wedge Y
  \enddiagram
  $$
Here the objects $X\wedge Y$ and $Y\wedge X$ are objects of the
first iteration, that is they are vertices from $V_1$.

$2^3$-diagram looks as follows:
  $$
  \diagram
  X\wedge X\wedge X \ar[dddddd]^-{1\wedge f\wedge 1}
  \ar[rrrr]^-{f\wedge 1\wedge 1}
  \ar[ddr]^-{1\wedge 1\wedge f} & & & &
  Y\wedge X\wedge X \ar[dddddd]^-{1\wedge f\wedge 1}
  \ar[ddr]^-{1\wedge 1\wedge f} & & \\ \\
  & X\wedge X\wedge Y \ar[dddddd]^-{1\wedge f\wedge 1}
  \ar[rrrr]^-{f\wedge 1\wedge 1}
  & & & & Y\wedge X\wedge Y \ar[dddddd]^-{1\wedge f\wedge 1} & \\ \\ \\ \\
  X\wedge Y\wedge X \ar[rrrr]^-{f\wedge 1\wedge 1}
  \ar[ddr]^-{1\wedge 1\wedge f} & & & & Y\wedge Y\wedge X
  \ar[ddr]^-{1\wedge 1\wedge f} & & \\ \\
  & X\wedge Y\wedge Y \ar[rrrr]^-{f\wedge 1\wedge 1}
  & & & & Y\wedge Y\wedge Y &
  \enddiagram
  $$
Here the objects $X\wedge Y\wedge Y$, $Y\wedge X\wedge Y$ and
$Y\wedge Y\wedge X$ are objects of the first iteration and the
objects $X\wedge Y\wedge X$, $X\wedge X\wedge Y$ and $Y\wedge
X\wedge X$ are objects of the second iteration.

\subsection{\it Mixed idempotents and their images}
\label{mixed}

Now let $\sg _n$ be a symmetric group of permutations of an
$n$-set and let $\sigma $ be any permutation from $\sg _n$. Let
$i\in \{ 0,1,\dots ,m\} $ and let $S\in \mathcal S_i$. Let $A_S=\{
a_1,\dots ,a_m\} $ be an ordered binary $m$-set, such that, for
any $j\in \{ 1,\dots ,m\} $, $a_j=0$ if $j$-th factor in the
product $(Y,X)_S$ coincides with $X$, and $a_j=1$ otherwise, i.e.
when its coincides with $Y$. Then the permutation $\sigma $
induces a uniquely defined permutation $\sigma _S$ on $A_S$. We
will say that the product $(Y,X)_S$ is of an inner type with
respect to the permutation $\sigma $ if $\sigma _S(a_j)=a_j$ for
any $j\in \{ 1,\dots ,m\} $, and we will say that $(Y,X)_S$ is of
an outer type with respect to $\sigma $ if there exists $j\in \{
1,\dots ,m\} $, such that $\sigma _S(a_j)\neq a_j$.

If $(Y,X)_S$ is of the inner type with respect to $\sigma $, then
let
  $$
  \Gamma _{\sigma ,S}:(Y,X)_S\lra (Y,X)_S
  $$
be a uniquely defined automorphism induced by $\sigma $ using
commutativity and associativity axioms in the monoidal category
$\C $. We may say again that $\Gamma _{\sigma ,S}$ is a graph of
the permutation $\sigma $. If $(Y,X)_S$ is of the outer type with
respect to $\sigma $, then we still can define a graph of $\sigma
$ as a uniquely defined isomorphism
  $$
  \Gamma _{\sigma ,S}:(Y,X)_S\lra (Y,X)_{\sigma (S)}\; ,
  $$
again induced by $\sigma $ using commutativity and associativity
in $\C $.

Let
  $$
  \Gamma _{\sigma ,i}':
  \coprod _{S\in \mathcal S_i}
  (Y,X)_S\lra
  \coprod _{S\in \mathcal S_i}
  (Y,X)_S
  $$
a morphism on coproducts induced by all maps $\Gamma _{\sigma ,S}$
when $S$ runs $\mathcal S_i$, but $\sigma $ is fixed. Then, for
any $i\in \{ 1,\dots ,m\} $, we have the commutative diagram
  $$
  \diagram
  \coprod _{S\in \mathcal S_i}(Y,X)_S
  \ar[rr]^-{} \ar[dd]^-{\Gamma _{\sigma ,i}'} & &
  \coprod _{S\in \mathcal S_{i-1}}(Y,X)_S
  \ar[dd]^-{\Gamma _{\sigma ,i-1}'} \\ \\
  \coprod _{S\in \mathcal S_i}(Y,X)_S
  \ar[rr]^-{} & &
  \coprod _{S\in \mathcal S_{i-1}}(Y,X)_S
  \enddiagram
  $$
where horizontal arrows are induced by all possible maps from the
objects of $V_i$ into the objects of $V_{i-1}$. In particular, we
have the commutative diagram
  $$
  \diagram
  X^{(m)} \ar[rr]^-{t_{m,i}} \ar[dd]^-{\Gamma _{\sigma ,m}'} & &
  \coprod _{S\in \mathcal S_i}(Y,X)_S \ar[dd]^-{\Gamma _{\sigma ,i}'} \\ \\
  X^{(m)} \ar[rr]^-{t_{m,i}} & & \coprod _{S\in \mathcal S_i}(Y,X)_S
  \enddiagram
  $$
Therefore, $\Gamma _{\sigma ,i}'$ induces an endomorphism
  $$
  \Gamma _{\sigma ,i}:(Y,X,m-i,i)\lra (Y,X,m-i,i)
  $$
on the level of colimits and, moreover, we have the commutative
diagram
  $$
  \diagram
  (Y,X,m-i,i)\ar[rr]^-{w_{m,i}} \ar[dd]^-{\Gamma _{\sigma ,i}} & &
  (Y,X,m-i+1,i-1) \ar[dd]^-{\Gamma _{\sigma ,i-1}} \\ \\
  (Y,X,m-i,i)\ar[rr]^-{w_{m,i}} & &
  (Y,X,m-i+1,i-1)
  \enddiagram
  $$
for each $i\in \{ 1,\dots ,m\} $.

In the same fashion, any product $(Z,X)_S$ can be of inner or of
outer type with respect to a permutation $\sigma \in \sg _m$. If
$(Z,X)_S$ is of the inner type, then let
  $$
  \Xi _{\sigma ,S}:(Z,X)_S\lra (Z,X)_S
  $$
be a unique automorphism induced by $\sigma $, and if $(Y,X)_S$ is
of the outer type, then we can define its graph as a uniquely
defined isomorphism
  $$
  \Xi _{\sigma ,S}:(Z,X)_S\lra (Z,X)_{\sigma (S)}\; ,
  $$
induced by $\sigma $ and the commutativity and associativity in
$\C $. Let
  $$
  \Xi _{\sigma ,i}':
  \vee _{S\in \mathcal S_i}
  (Z,X)_S\lra
  \vee _{S\in \mathcal S_i}
  (Z,X)_S
  $$
a morphism on coproducts induced by all maps $\Xi _{\sigma ,S}$
when $S\in \mathcal S_i$ and let
  $$
  \Xi _{\sigma ,i}:[Z,X,m-i,i]\lra [Z,X,m-i,i]
  $$
be the corresponding morphism in the category $\T $. Then, for any
$i\in \{ 1,\dots ,m\} $ we have the morphism of cofibered
sequences
  $$
  \diagram
  (Y,X,m-i,i) \ar[r]^-{w_{m,i}} \ar[dd]^-{\Gamma _{\sigma ,i}} &
  (Y,X,m-i+1,i-1) \ar[r]^-{} \ar[dd]^-{\Gamma _{\sigma ,i-1}} &
  [Z,X,m-i+1,i-1] \ar[dd]^-{\Xi _{\sigma ,i-1}}  \\ \\
  (Y,X,m-i,i) \ar[r]^-{w_{m,i}} &
  (Y,X,m-i+1,i-1) \ar[r]^-{} &
  [Z,X,m-i+1,i-1]
  \enddiagram
  $$
Applying Lemma \ref{natur}, we obtain the corresponding morphisms
of distinguished triangles

{\footnotesize

  $$
  \diagram
  (Y,X,m-i,i) \ar[r]^-{w_{m,i}} \ar[dd]^-{\Gamma _{\sigma ,i}} &
  (Y,X,m-i+1,i-1) \ar[r]^-{} \ar[dd]^-{\Gamma _{\sigma ,i-1}} &
  [Z,X,m-i+1,i-1] \ar[r]^-{} \ar[dd]^-{\Xi _{\sigma ,i-1}} &
  \sg (Y,X,m-i,i) \ar[dd]^-{\sg \Gamma _{\sigma ,i}} \\ \\
  (Y,X,m-i,i) \ar[r]^-{w_{m,i}} &
  (Y,X,m-i+1,i-1) \ar[r]^-{} &
  [Z,X,m-i+1,i-1] \ar[r]^-{} &
  \sg (Y,X,m-i,i)
  \enddiagram
  $$

}

Let
  $$
  d^+_{m,i}=\frac{1}{m!}\sum _{\sigma \in \sg _m}
  sgn(\sigma )\Gamma _{\sigma ,i}\; ,
  $$
  $$
  d^-_{m,i}=\frac{1}{m!}\sum _{\sigma \in \sg _m}
  \Gamma _{\sigma ,i}\; ,
  $$
  $$
  e^+_{m,i}=\frac{1}{m!}\sum _{\sigma \in \sg _m}
  sgn(\sigma )\Xi _{\sigma ,i}
  $$
and
  $$
  e^-_{m,i}=\frac{1}{m!}\sum _{\sigma \in \sg _m}
  \Xi _{\sigma ,i}
  $$
for any $i\in \{ 0,1,\dots ,m\} $. It is not hard to see that all
of these maps are idempotents in $\T $. Note that $d^{\pm
}_{m,0}=d^{\pm }_m$ for $Y$, and $d^{\pm }_{m,m}=d^{\pm }_m$ for
$X$, where $d^{\pm }_m$ are the idempotents defined in Section
\ref{basics}. Similarly, $e^{\pm }_{m,0}=d^{\pm }_m$ for $Z$, and
$e^{\pm }_{m,m}=d^{\pm }_m$ for $X$. Therefore we will say that
$d^{\pm }_m$ are pure idempotents and that $d^{\pm }_{m,i}$ and
$e^{\pm }_{m,i}$ are mixed idempotents.

Now, summing vertical maps in the last morphism of distinguished
triangles, we obtain two morphisms of distinguished triangles

{\footnotesize

  $$
  \diagram
  (Y,X,m-i,i) \ar[r]^-{w_{m,i}} \ar[dd]^-{d^+_{m,i}} &
  (Y,X,m-i+1,i-1) \ar[r]^-{} \ar[dd]^-{d^+_{m,i-1}} &
  [Z,X,m-i+1,i-1] \ar[r]^-{} \ar[dd]^-{e^+_{m,i-1}} &
  \sg (Y,X,m-i,i) \ar[dd]^-{\sg d^+_{m,i}} \\ \\
  (Y,X,m-i,i) \ar[r]^-{w_{m,i}} &
  (Y,X,m-i+1,i-1) \ar[r]^-{} &
  [Z,X,m-i+1,i-1] \ar[r]^-{} &
  \sg (Y,X,m-i,i)
  \enddiagram
  $$

}

\noindent and

{\footnotesize

  $$
  \diagram
  (Y,X,m-i,i) \ar[r]^-{w_{m,i}} \ar[dd]^-{d^-_{m,i}} &
  (Y,X,m-i+1,i-1) \ar[r]^-{} \ar[dd]^-{d^-_{m,i-1}} &
  [Z,X,m-i+1,i-1] \ar[r]^-{} \ar[dd]^-{e^-_{m,i-1}} &
  \sg (Y,X,m-i,i) \ar[dd]^-{\sg d^-_{m,i}} \\ \\
  (Y,X,m-i,i) \ar[r]^-{w_{m,i}} &
  (Y,X,m-i+1,i-1) \ar[r]^-{} &
  [Z,X,m-i+1,i-1] \ar[r]^-{} &
  \sg (Y,X,m-i,i)
  \enddiagram
  $$

}

To extract images of these idempotents of distinguished triangles
we need some trivial, but useful lemmas.

Let $\X $ be a category and let $f:X\to X$ be an idempotent in $\X
$, i.e. $f^2=f$. The typical example: if $g:X\to Y$ and $h:Y\to X$
are two morphisms in $\X $, such that $gh=1_Y$, then $f=hg$ is an
idempotent. In such a case $f$ is called to be a splitting
idempotent. If the idempotent $f$ has two splittings $f=hg=st$.
Then there exists a unique isomorphism $d:Y\stackrel{\cong }{\to
}Z$, such that the diagram
  $$
  \diagram
  & & Y\ar[dd]^-{d} \ar[rrd]^-{h} & & \\
  X\ar[rru]^-{g} \ar[rrd]^-{t} & & & & X \\
  & & Z \ar[rru]^-{s} & &
  \enddiagram
  $$
commutes. Indeed, define $d$ as $th$ and set $c=gs$. Then
$dc=thgs=tfs=tsts=1_Z$ and $cd=gsth=gsth=ghgh=1_Y$. Furthermore,
$dg=thg=tf=tst=t$ and $sd=sth=fh=hgh=h$. If there exists another
such $d'$, then $sd=h=sd'$, whence $tsd=tsd'$. Then $d=d'$ because
of $ts=1_Z$. So, we may speak about the image $im(f)$ of a
splitting idempotent $f$ which is defined up to a canonical
isomorphism.

In fact, this is a corollary of more general

\begin{lemma}
\label{imagefunctor} Let $f:X\to X$ and $g:Y\to Y$ be splitting
idempotents in $\X $ and let $t:X\to Y$ be a morphism from $f$ to
$g$, i.e. $gt=tf$. Then there exists a unique morphism $q:im(f)\to
im(g)$, such that the diagram
  $$
  \diagram
  X \ar[rr]^-{a} \ar[dd]^-{t} & & im(f) \ar[dd]^-{q}
  \ar[rr]^-{b} & & X \ar[dd]^-{t} \\ \\
  Y \ar[rr]^-{c} & & im(g) \ar[rr]^-{d} & & Y
  \enddiagram
  $$
commutes.
\end{lemma}

\begin{pf}
To see the uniqueness assume that $q$ exists. Then $qa=ct$. Since
$ab=1_X$, it follows that $q$ is uniquely determined as $q=ctb$.
Defining $q$ as $ctb$ we get: $qa=ctba=ctf=cgt=cdct=ct$ and
$dq=dctb=gtb=tfb=tbab=tb$.
\end{pf}

One can say that an image of idempotents is a functor on the
evident category of splitting idempotents over the category $\X $.

\begin{lemma}
\label{triangleidemp} Let
  $$
  \diagram
  X \ar[rr]^-{t} \ar[dd]^-{f} & & Y \ar[rr]^-{s} \ar[dd]^-{g}
  & & Z \ar[rr]^-{u} & & \sg X \ar[dd]^-{\sg f} \\ \\
  X \ar[rr]^-{t} & & Y \ar[rr]^-{s} & & Z \ar[rr]^-{u} & & \sg X
  \enddiagram
  $$
be a commutative diagram, where $f$ and $g$ are splitting
idempotents and both rows are the same distinguished triangle.
Then there exists a splitting idempotent $h:Z\to Z$ (not necessary
unique), such that the diagram
  $$
  \diagram
  X \ar[rr]^-{t} \ar[dd]^-{} & & Y \ar[rr]^-{s} \ar[dd]^-{}
  & & Z \ar[rr]^-{u} \ar[dd]^-{} & & \sg X \ar[dd]^-{} \\ \\
  im(f) \ar[rr]^-{} \ar[dd]^-{} & & im(g) \ar[rr]^-{} \ar[dd]^-{}
  & & im(h) \ar[rr]^-{} \ar[dd]^-{} & & \sg (im(f)) \ar[dd]^-{} \\ \\
  X \ar[rr]^-{t} & & Y \ar[rr]^-{s} & & Z
  \ar[rr]^-{u} & & \sg X
  \enddiagram
  $$
is commutative, where all vertical compositions are splittings of
the corresponding idempotents. Moreover, the middle triangle is
distinguished.
\end{lemma}

\begin{pf}
Indeed, let $C$ be a cone of the morphism $im(t):im(f)\to im(g)$
and let $a$ and $b$ be two morphisms on cones included in the
commutative diagram
  $$
  \diagram
  im(f) \ar[rr]^-{} \ar[dd]^-{} & & im(g) \ar[rr]^-{} \ar[dd]^-{}
  & & C \ar[rr]^-{} \ar[dd]^-{a} & & \sg (im(f)) \ar[dd]^-{} \\ \\
  X \ar[rr]^-{} \ar[dd]^-{} & & Y \ar[rr]^-{} \ar[dd]^-{}
  & & Z \ar[rr]^-{} \ar[dd]^-{b} & & \sg X \ar[dd]^-{} \\ \\
  im(f) \ar[rr]^-{} & & im(g) \ar[rr]^-{s} & & C
  \ar[rr]^-{} & & \sg (im(f))
  \enddiagram
  $$
Since the compositions $im(f)\to X\to im(f)$, $im(g)\to Y\to
im(g)$ are the identity, $c=ba$ is an isomorphism. Let
$d=c^{-1}b$. Then we have the commutative diagram
  $$
  \diagram
  im(f) \ar[rr]^-{} \ar[dd]^-{} & & im(g) \ar[rr]^-{} \ar[dd]^-{}
  & & C \ar[rr]^-{} \ar[dd]^-{a} & & \sg (im(f)) \ar[dd]^-{} \\ \\
  X \ar[rr]^-{} \ar[dd]^-{} & & Y \ar[rr]^-{} \ar[dd]^-{}
  & & Z \ar[rr]^-{} \ar[dd]^-{d} & & \sg X \ar[dd]^-{} \\ \\
  im(f) \ar[rr]^-{} & & im(g) \ar[rr]^-{s} & & C
  \ar[rr]^-{} & & \sg im(f)
  \enddiagram
  $$
where all three vertical compostions are identical morphisms. Then
$h=ad$ and $C=im(h)$.

The triangle
  $$
  im(f)\longrightarrow im(g)\longrightarrow
  im(h)\longrightarrow \sg im(f)
  $$
is a candidate triangle in $\T $, see \cite{Nee}, i.e. all three
compositions in it are trivial. This is so because, if we consider
the commutative diagram
  $$
  \diagram
  X \ar[rr]^-{t} \ar[dd]^-{f'} & & Y \ar[rr]^-{s} \ar[dd]^-{g'}
  & & Z \ar[rr]^-{u} \ar[dd]^-{h'} & & \sg X \ar[dd]^-{\sg f'} \\ \\
  im(f) \ar[rr]^-{im(t)} & & im(g) \ar[rr]^-{im(s)} & & im(h)
  \ar[rr]^-{im(u)} & & \sg im(f)
  \enddiagram
  $$
then, say, $im(s)\circ im(t)\circ f'=h'\circ s\circ t=0$, whence
$im(s)\circ im(t)=0$ because $f'$ has a section.

By symmetry the same holds for the triangle
  $$
  im(1_X-f)\longrightarrow im(1_Y-g)\longrightarrow
  im(1_Z-h)\longrightarrow \sg im(1_X-f)
  $$
At the same time, evidently, the dist. triangle $X\to Y\to Z\to
\sg X$ is a direct sum of the two above candidate triangles.
Therefore, both candidate triangles are distinguished, see
\cite{Nee}.
\end{pf}

\begin{proposition}
\label{triangleimages} There are commutative diagrams

{\footnotesize

  $$
  \diagram
  (Y,X,m-i,i) \ar[r]^-{w_{m,i}} \ar[dd]^-{} & (Y,X,m-i+1,i-1)
  \ar[r]^-{} \ar[dd]^-{}
  & [Z,X,m-i+1,i-1] \ar[r]^-{} \ar[dd]^-{} & \sg (Y,X,m-i,i) \ar[dd]^-{} \\ \\
  I^+_{m,i} \ar[r]^-{w^+_{m,i}} \ar[dd]^-{} & I^+_{m,i-1}
  \ar[r]^-{} \ar[dd]^-{}
  & J^+_{m,i-1} \ar[r]^-{} \ar[dd]^-{} & \sg I^+_{m,i} \ar[dd]^-{} \\ \\
  (Y,X,m-i,i) \ar[r]^-{w_{m,i}} & (Y,X,m-i+1,i-1) \ar[r]^-{} &
  [Z,X,m-i+1,i-1]
  \ar[r]^-{} & \sg (Y,X,m-i,i)
  \enddiagram
  $$

}

\noindent and

{\footnotesize

  $$
  \diagram
  (Y,X,m-i,i) \ar[r]^-{w_{m,i}} \ar[dd]^-{} & (Y,X,m-i+1,i-1)
  \ar[r]^-{} \ar[dd]^-{}
  & [Z,X,m-i+1,i-1] \ar[r]^-{} \ar[dd]^-{} & \sg (Y,X,m-i,i) \ar[dd]^-{} \\ \\
  I^-_{m,i} \ar[r]^-{w^-_{m,i}} \ar[dd]^-{} & I^-_{m,i-1}
  \ar[r]^-{} \ar[dd]^-{}
  & J^-_{m,i-1} \ar[r]^-{} \ar[dd]^-{} & \sg I^-_{m,i} \ar[dd]^-{} \\ \\
  (Y,X,m-i,i) \ar[r]^-{w_{m,i}} & (Y,X,m-i+1,i-1) \ar[r]^-{} &
  [Z,X,m-i+1,i-1]
  \ar[r]^-{} & \sg (Y,X,m-i,i)
  \enddiagram
  $$

}

\noindent where
 $$
 I^{\pm }_{m,i}=im(d^{\pm }_{m,i})
 $$
and
 $$
 J^{\pm }_{m,i}=im(e^{\pm }_{m,i})
 $$
are images of the mixed idempotents and the maps $w^+_{m,i}$ and
$w^-_{m,i}$ are induced on images of idempotents by the map
$w_{m,i}$. All rows in these commutative diagrams are
distinguished triangles and all columns are splittings of the
corresponding idempotents $d^{\pm }_{m,i}$ and $e^{\pm }_{m,i}$.
\end{proposition}

\begin{pf}
Apply Lemma \ref{triangleidemp} to the above maps of distinguished
triangles induced by the mixed idempotents $d^{\pm }_{m,i}$ and
$e^{\pm }_{m,i}$.
\end{pf}

Now we want to compute the images $J^{\pm }_{m,i}$ of idempotents
$e^{\pm }_{m,i}$. We again need an abstract lemma.

\begin{lemma}
\label{isoimage} Let $\X $ be a $\mathbb Q$-linear category with
splitting idempotents, and let
  $$
  \diagram
  A \ar[rr]^-{u} \ar[dd]^-{a} & & B \ar[dd]^-{b} \\ \\
  A & & B \ar[ll]^-{d}
  \enddiagram
  $$
be a diagram in $\X $, commutative up to a scalar $\alpha $, i.e.
$\alpha a=dbu$ for some $\alpha \in \mathbb Q$. Assume,
furthermore, that $b^2=b$ and $ud=\alpha $ (to be more precise,
$ud=\alpha 1_B$). Then $a^2=a$, $ua=bu$, $ad=db$ and the map
$im(u):im(a)\to im(b)$ (which exists because $ua=bu$, see Lemma
\ref{imagefunctor}) is an isomorphism.
\end{lemma}

\begin{pf}
Indeed, since $\alpha a=dbu$ and $ud=\alpha $, it follows that
$\alpha ua=\alpha bu$, whence $ua=bu$. Similarly, multiplying
$\alpha a=dbu$ on $d$ from the right hand side, we have that
$\alpha ad=dbud$. Since $ud=\alpha $, we get $\alpha ad=db\alpha
$, whence $ad=db$.

Further, $a^2=\alpha ^{-2}(dbu)(dbu)=\alpha ^{-2}db(ud)bu=\alpha
^{-2}db\alpha bu=\alpha ^{-1}dbbu$. Since $b^2=b$ by assumtions,
we get $a^2=\alpha ^{-1}dbu=a$.

Now let's consider the commutative diagram
  $$
  \diagram
  B \ar[rr]^-{d} \ar[dd]^-{\pi _B} & &
  A \ar[rr]^-{u} \ar[dd]^-{\pi _A} & &
  B \ar[rr]^-{d} \ar[dd]^-{\pi _B} & &
  A \ar[dd]^-{\pi _A} \\ \\
  I_B \ar[rr]^-{I_d} \ar[dd]^-{\iota _B} & &
  I_A \ar[rr]^-{I_u} \ar[dd]^-{\iota _A} & &
  I_B \ar[rr]^-{I_d} \ar[dd]^-{\iota _B} & &
  I_A \ar[dd]^-{\iota _A} \\ \\
  B \ar[rr]^-{d}  & &
  A \ar[rr]^-{u}  & &
  B \ar[rr]^-{d}  & &
  A
  \enddiagram
  $$
where the columns are splittings of the idempotents $a$ and $b$,
$I_A=im(a)$, $I_B=im(b)$, etc. Easy chasing on this diagram shows
that $\iota _BI_uI_d\pi _B=bud=\alpha \iota _B\pi _B$, whence
$\iota _B(\alpha ^{-1}I_uI_d)\pi _B=\iota _B\pi _B$. Since $\iota
_B$ is a left inverse for $\pi _B$, it follows that $\alpha
^{-1}I_uI_d=1_{I_B}$. Further, $I_dI_u=(\pi _Ad\iota _B)(\pi
_Bu\iota _A)=\pi _Ad(\iota _B\pi _B)u\iota _A=\pi _Adbu\iota
_A=\pi _A(\alpha a)\iota _A=\alpha \pi _Aa\iota _A=\alpha \pi
_A\iota _A\pi _A\iota _A=\alpha $, whence $\alpha
^{-1}I_dI_u=1_{I_A}$.
\end{pf}

\begin{proposition}
\label{imcompute}
$$
J^+_{m,i}\cong Z^{[m-i)}\otimes X^{[i)}
$$
$$
J^-_{m,i}\cong Z^{(m-i]}\otimes X^{(i]}
$$
\end{proposition}

\begin{pf}
Let $\sg _m$ be a symmetric group of permutations of an $m$-set
$\{ 1,\dots ,m\} $. Let $i$ be an integer, such that $1\leq i\leq
m$. The group $\sg _i$ may be considered as a subgroup in $\sg _m$
consisting of permutations which preserve the set $\{ 1,\dots ,i\}
$, and $\sg _{m-i}$ as a subgroup in $\sg _m$ of permutations
preserving the set $\{ i+1,\dots ,m\} $. Then $\sg _{m-i}\times
\sg _i$ is a subgroup in $\sg _m$ consisting of the product of
permutations acting on the sets $\{ 1,\dots ,i\} $ and $\{
i+1,\dots ,m\} $ in an inner fashion.

The number of elements in $\mathcal S_i$ is equal to the number
$C_m^r=\frac{m!}{(m-i)!i!}$, so that it coincides with the number
of left cosets in $\sg _m$ modulo the subgroup $\sg _{m-i}\times
\sg _i$. For any $S\in \mathcal S_i$ fix an isomorphism
  $$
  u_S:[Z,X]_S\stackrel{\cong }{\lra }
  Z^{(m-i)}\otimes X^{(i)}
  $$
fixing, in fact, a representative $\varsigma _S$ of the correspond
ing right coset in $\sg _m$ modulo the subgroup $\sg _{m-i}\times
\sg _i$. Let
  $$
  u:[Z,X,m-i,i]\stackrel{}{\lra }
  Z^{(m-i)}\otimes X^{(i)}
  $$
be a sum of isomorphisms $u_S$, $S\in \mathcal S_i$, and let
  $$
  d:Z^{(m-i)}\otimes X^{(i)}
  \stackrel{}{\lra }[Z,X,m-i,i]
  $$
be a morphism, such that its composition with any projection
  $$
  p_S:[Z,X,m-i,i]\lra [Z,X]_S
  $$
coincides with the corresponding inverse isomorphism $u_S^{-1}$.
Let's consider the diagram
  $$
  \diagram
  [Z,X,m-i,i] \ar[rr]^-{u} \ar[dd]^-{e^{\pm }_{m,i}} & &
  Z^{(m-i)}\otimes X^{(i)}
  \ar[dd]^-{d^{\pm }_{m-i}\otimes d^{\pm }_i} \\ \\
  [Z,X,m-i,i] & &
  Z^{(m-i)}\otimes X^{(i)}
  \ar[ll]^-{d}
  \enddiagram
  $$
Since $\{ \varsigma _S\} _{S\in \mathcal S_i}$ is a set of
representatives of the right cosets in $\sg _m$ modulo the
subgroup $\sg _{m-i}\times \sg _i$, $\{ \varsigma _S^{-1}\} _{S\in
\mathcal S_i}$ is a collection of representatives of the right
cosets $\sg _m$ modulo $\sg _{m-i}\times \sg _i$. Then
  $$
  \sg _m=\sg _m\varsigma _{S}=\cup _{T\in \mathcal S_i}
  \varsigma _T^{-1}(\sg _{m-i}\times \sg _i)\varsigma _S\; ,
  $$
so that
  $$
  (m-i)!i!(p_T\circ d\circ (d^{\pm }_{m-i}\otimes d^{\pm }_i)\circ
  u_S)=m!(p_T\circ e^{\pm }_{m,i}\circ u_S)\; ,
  $$
or, equivalently,
  $$
  p_T\circ d\circ (d^{\pm }_{m-i}\otimes d^{\pm }_i)\circ
  u_S=\frac{m!}{(m-i)!i!}p_T\circ e^{\pm }_{m,i}\circ u_S
  $$
for any $S$ and $T$ from $\mathcal S_i$. This shows that the above
diagram is commutative modulo the scalar $\frac{m!}{(m-i)!i!}$,
i.e.
  $$
  \frac{m!}{(m-i)!i!}\cdot e^{\pm }_{m,i}=
  d\circ (d^{\pm }_{m-i}\otimes d^{\pm }_i)\circ u\; .
  $$
Moreover, the composition $ud$, evidently, coincides with the
multiplication by $\frac{m!}{(m-i)!i!}$. Now it remains just to
apply Lemma \ref{isoimage} and observe that
  $$
  im(d^{\pm }_{m-i}\otimes d^{\pm }_i)=
  im(d^{\pm }_{m-i})\otimes im(d^{\pm }_i)\; .
  $$
\end{pf}

\subsection{\it Triangle filtrations: the proof of Theorem \ref{main1}}

Let's recall that the following approach to Theorem \ref{main1}
was suggested by Uwe Jannsen: build a filtration on wedge
(symmetric) powers of vertices in a distinguished triangle,
similar to the filtration for powers in a short exact sequence of
locally free sheaves of modules on a manifold. Now it remains just
to claim that the middle rows in the diagrams from Proposition
\ref{triangleimages} provide, in fact, the desired filtration.

Let $\T $ be an arbitrary triangulated category and let $X$ be an
object in $\T $. The following definition appears in \cite{Ka}, p.
152 - 153. A finite (decreasing) filtration on $X$ is a sequence
of objects $A_{-1},A_0,A_1,A_2,\dots $ and of distinguished
triangles
  $$
  F^iX\; : \; \; A_i\stackrel{a_i}{\lra }A_{i-1}\lra Gr_F^{i-1}X\lra \sg A_i\; ,
  $$
such that
  $$
  F^0X\; : \; \; X\stackrel{1}{\lra }X\lra 0\lra \sg X\; ,
  $$
i.e. $A_{-1}=X$, $A_0=X$ and $a_0=1_X$, and
  $$
  F^{m+1}X\; : \; \; 0\lra A_m\stackrel{1}{\lra }A_{m-1}\lra 0\; ,
  $$
for some natural number $m$, i.e. $A_{m+1}=0$ and $Gr_F^mX= A_m$.
The last condition may be interpreted also like $F^{m+1}=0$. Note
that if graded pieces $Gr_F^iX$ are equal to zero for all $i$,
then it follows that $X$ is also trivial.

\begin{proposition}
\label{filt} Let $\T =Ho(\C )$ be a trinagulated category,
reinforced by an underlying pointed model and monoidal category
$\C $ in the sense of Section \ref{reinforce}, and let
 $$
 X\stackrel{f}{\lra }Y\stackrel{g}{\lra }Z\stackrel{h}{\lra }\sg X
 $$
be a distinguished triangle in $\T $. Then, for any natural $m$,
there exist a finite decreasing filtration $F^*Y^{[m)}$ on
$Y^{[m)}$ and a finite decreasing filtration $F^*Y^{(m]}$ on
$Y^{(m]}$, such that
  $$
  Gr^i_FY^{[m)}\cong Z^{[m-i)}\otimes X^{[i)}
  $$
and
  $$
  Gr^i_FY^{(m]}\cong Z^{(m-i]}\otimes X^{(i]}
  $$
for any $i$, and, in both cases, $F^{m+1}=0$.
\end{proposition}

\begin{pf}
Recall that, without loss of generality, applying cofibrant
replacement, we may assume that $f$ is a cofibration and both $X$
and $Y$ are cofibrant, so that the above dist. triangle is a
cofibration triangle. Then we define a filtration on $Y^{[m)}$
using distinguished triangles of images of mixed idempotents from
Section \ref{mixed}. In other words, let
  $$
  F^0Y^{[m)}\; : \; \; Y^{[m)}\stackrel{1}{\lra }
  Y^{[m)}\lra 0\lra \sg Y^{[m)}\; ,
  $$
and let
  $$
  F^iY^{[m)}\; : \; \; I^+_{m,i}\stackrel{w^+_{m,i}}{\lra }
  I^+_{m,i-1}\lra J^+_{m,i-1}\lra \sg I^+_{m,i}
  $$
for all $i\in \{ 1,\dots ,m\} $. Also, let
  $$
  F^{m+1}Y^{[m)}\; : \; \; 0\lra
  I^+_{m,m}\stackrel{1}{\lra }I^+_{m,m}\lra 0
  $$
be the last term of the filtration. Then it is easy to see that
the first term of this filtration is
  $$
  F^1Y^{[m)}\; : \; \; I^+_{m,1}\stackrel{w^+_{m,1}}{\lra }
  Y^{[m)}\stackrel{g^{[m)}}{\lra }Z^{[m)}\lra \sg I^+_{m,1}
  $$
and the last one is the triangle
  $$
  F^{m+1}Y^{[m)}\; : \; \; 0\lra
  X^{[m)}\stackrel{1}{\lra }X^{[m)}\lra 0
  $$
By Prop. \ref{imcompute}, the graded pieces can be computed by the
formula
  $$
  Gr_F^iY^{[m)}=J^+_{m,i}\cong Z^{[m-i)}\otimes X^{[i)}
  $$
for any $i$.

Similarly, we define a filtration on $Y^{(m]}$:
  $$
  F^0Y^{(m]}\; : \; \; Y^{(m]}\stackrel{1}{\lra }
  Y^{(m]}\lra 0\lra \sg Y^{(m]}\; ,
  $$
  $$
  F^iY^{(m]}\; : \; \; I^-_{m,i}\stackrel{w^-_{m,i}}{\lra }
  I^-_{m,i-1}\lra J^-_{m,i-1}\lra \sg I^-_{m,i}
  $$
for all $i\in \{ 1,\dots ,m\} $, and set
  $$
  F^{m+1}Y^{(m]}\; : \; \; 0\lra
  I^-_{m,m}\stackrel{1}{\lra }I^-_{m,m}\lra 0
  $$
as a last term of the filtration. Then
  $$
  F^1Y^{(m]}\; : \; \; I^-_{m,1}\stackrel{w^-_{m,1}}{\lra }
  Y^{(m]}\stackrel{g^{(m]}}{\lra }Z^{(m]}\lra \sg I^-_{m,1}
  $$
and
  $$
  F^{m+1}Y^{(m]}\; : \; \; 0\lra
  X^{(m]}\stackrel{1}{\lra }X^{(m]}\lra 0
  $$
The graded pieces:
  $$
  Gr_F^iY^{(m]}=J^-_{m,i}\cong Z^{(m-i]}\otimes X^{(i]}
  $$
for any $i$.
\end{pf}

\bigskip

Now we can finish the proof of Theorem \ref{main1}. Assume that
$X$ and $Z$ in the triangle from Theorem \ref{main1} are evenly
finite dimensional. It means that $X^{[t)}=0$ and $Z^{[t)}=0$ for
some natural $t$. Then, for a big enough $m$, all graded pieces
$Gr_F^iY^{[m)}\cong Z^{[m-i)}\otimes X^{[i)}$ of the even
filtration are equal to zero. Therefore, $Y^{[m)}=0$. And
similarly in the odd case.

As to the the second part of Theorem \ref{main1}, it is equivalent
to the first one. To see this we need just to observe that the
shift functor $\sg $ in any $\mathbb Q$-linear monoidal
triangulated category $\T $ carries evenly (oddly) finite
dimensional objects into oddly (evenly) finite dimensional
objects. This is due to the axioms coding the compatibility of the
monoidal and the triangulated structures in $\T $, see A8 in
\cite{VMW}.

\section{Finite dimensionality of motives of curves over a field}

In this section we prove Theorem \ref{main2}. The word "scheme"
means a separated scheme of finite type over a field $k$ and the
word "curve" means a an integral one-dimensional scheme over $k$.
We will also assume that $char(k)=0$.

\subsection{\it Scalar extension and a splitting lemma}

If $k$ is algebraically closed, then $X$ can be considered as a
Zariski open subset in a projective curve $Y$, see \cite{Ha}, p.
105, and \cite{Nag}. Let $p:W\to Y$ be a resolution of
singularities of $Y$ and let $U=p^{-1}(X)$, so that we have a
commutative square
  \begin{equation}
  \label{sing}
  \diagram
  W \ar[rr]^-{p} & & Y \\ \\
  U \ar[rr]^-{p} \ar[uu]^-{} & & X \ar[uu]^-{}
  \enddiagram
  \end{equation}
Let also $Z=Y-X$ and $V=W-U$ be complements of Zariski open
subsets $Y$ and $U$ in the projective curves $X$ and $W$
respectively.

If $k$ is not algebraically closed, then we may consider the
square (\ref{sing}) first over an algebraic closure of $k$ and
then take a finite extension $L/k$, such that all varieties and
maps in (\ref{sing}), as well as $Z$ and $V$, are defined over
$L$, and $L$ contains $\sqrt {-1}$. But, since we work with
motives with coefficients in $\mathbb Q$, one can use transfer
arguments to show that finite dimensionality of the motive
$M(X_L)$ in $\DM (L)_{\mathbb Q}$ implies finite dimensionality of
$M(X)$ in $\DM (k)_{\mathbb Q}$. It means that, proving Theorem
\ref{main2}, we may assume, without loss of generality, that all
the data in the square (\ref{sing}) is rational over $k$.

We will also need the following useful

\begin{lemma}
\label{usefullemma} Let $\T $ be a triangulated category with the
shift functor $\sg $. Assume that we have a distinguished triangle
  $$
  \diagram
  A\oplus B \ar[rr]^-{
  \tiny
  \left(
  \begin{array}{cc}
  a & b \\
  c & d
  \end{array}
  \right)
  } & & A\oplus C
  \ar[r]^-{} & D \ar[r]^-{} &
  \sg (A\oplus B)
  \enddiagram
  $$
where $a$ is an automorphism of the object $A$. Then this triangle
is isomorphic to the direct sum of triangles
  $$
  A\stackrel{1}{\lra }A\lra 0 \lra \sg A
  $$
and
  $$
  B\stackrel{t}{\lra }C\lra D \lra \sg B
  $$
where $t=d-ca^{-1}b$.
\end{lemma}

\begin{pf}
This is just a reformulation of Lemma 1.2.4 in \cite{Nee}.
\end{pf}

\subsection{\it The proof of Theorem \ref{main2}}
\label{pfmain2}

We will consider three cases separately: (a) when a curve is
projective, but not necessary smooth; (b) when it is not
projective, but smooth, and (c) when it is not projective and,
probably, not smooth. Then, to prove Theorem \ref{main2}, we will
have just to join all together.

\begin{proposition}
\label{a} Assume that $X$ is projective. Then $M(X)$ is finite
dimensional and, moreover,
  $$
  M(X)\cong \mathbb Q\oplus G\oplus \mathbb Q(1)[2]\; ,
  $$
where $G$ is an oddly finite dimensional object in $\DM _{\mathbb
Q}$.
\end{proposition}

\begin{pf}
If $X$ is projective and smooth, then the proposition holds by Th.
\ref{findimcurves}. Assume that $X$ is singular. For simplicity,
we will consider the case when $X$ has only one singular point
(the other case can be proved by the same methods, but with more
cumbrous formulas). Then $p:U\to X$ contracts points $\{ u_1,\dots
,u_n\} $ onto a singular point in $X$. Let
  \begin{equation}
  \label{1p}
  {\mathbb Q^{\oplus }}^n\lra \mathbb Q\oplus M(U)
  \lra M(X)\lra {\mathbb Q^{\oplus }}^n[1]
  \end{equation}
be a blow up distinguished triangle corresponding to the map $p$,
see \cite{SV}, Th. 5.2. Note that ${\mathbb Q^{\oplus }}^n$ is
just a motive of the finite set $\{ u_1,\dots ,u_n\} $.

For any $i$ the composition
  $$
  \mathbb Q\to {\mathbb Q^{\oplus }}^n \lra
  \mathbb Q\oplus M(U)\to \mathbb Q\; ,
  $$
corresponding to the point $u_i$, is an isomorphism. By Lemma
\ref{usefullemma} the triangle (\ref{1p}) is isomorphic to a
direct sum of two distinguished triangles
  \begin{equation}
  \label{2p}
  {\mathbb Q^{\oplus }}^{n-1}\stackrel{t}{\lra }M(U)\lra
  M(X)\lra {\mathbb Q[1]^{\oplus }}^{n-1}
  \end{equation}
and
  $$
  \mathbb Q\lra \mathbb Q\lra 0\lra \mathbb Q[1]\; ,
  $$
In other words, we split the isomorphism $\mathbb Q\stackrel{\cong
}{\to }\mathbb Q$ corresponding to a point from $\{ u_1,\dots
,u_n\} $, say, to $u_1$.

Since $U$ is smooth projective, we have a decomposition
  $$
  M(U)=\mathbb Q\oplus M^1(U)\oplus \mathbb Q(1)[2]
  $$
induced by some $k$-rational point on $U$. Then we rewrite
(\ref{2p}) as follows:
  \begin{equation}
  \label{3p}
  {\mathbb Q^{\oplus }}^{n-1}\stackrel{t}{\lra }
  \mathbb Q\oplus M^1(U)\oplus \mathbb Q(1)[2]\lra
  M(X)\lra {\mathbb Q[1]^{\oplus }}^{n-1}
  \end{equation}

For any $i$, let $\nu _i:Spec(k)\lra U$ be the map corresponding
to the point $u_i$. Let also $\gamma :U\to Spec (k)$ be the
structure map for $U$. If $i>1$, then the composition
  $$
  \mathbb Q\to {\mathbb Q^{\oplus }}^{n-1}
  \stackrel{t}{\lra }M(U)\; ,
  $$
corresponding to the point $u_i$, coincides with the difference
$M(\nu _i)-M(\nu _1)$ (here we use the general expression for the
morphism $t$ given by Lemma \ref{usefullemma}). The projection
  $$
  M(U)=\mathbb Q\oplus M^1(U)\oplus \mathbb Q(1)[2]
  \lra \mathbb Q
  $$
is, in fact, the morphism $M(\gamma ):M(U)\to M(Spec(k))$.
Therefore, for any $u_i$, $i>1$, the composition
  $$
  \mathbb Q\to {\mathbb Q^{\oplus }}^{n-1}
  \stackrel{t}{\lra }M(U)\lra \mathbb Q
  $$
is equal to the difference $M(\gamma \nu _i)-M(\gamma \nu _1)$,
which is equal to zero. And, in addition, any map from $\mathbb Q$
to $\mathbb Q(1)[2]$ is zero. This shows that the triangle
(\ref{3p}) is a direct sum of two distinguished triangles
  \begin{equation}
  \label{4p}
  {\mathbb Q^{\oplus }}^{n-1}\stackrel{t}{\lra }M^1(U)\lra
  G\lra {\mathbb Q[1]^{\oplus }}^{n-1}
  \end{equation}
and
  $$
  0\lra \mathbb Q\oplus \mathbb Q(1)[2]\lra
  \mathbb Q\oplus \mathbb Q(1)[2]\lra 0\; .
  $$
In particular,
  $$
  M(X)=\mathbb Q\oplus G\oplus \mathbb Q(1)[2]\; .
  $$

Now recall that $M^1(U)$ is oddly finite dimensional by Th.
\ref{findimcurves}. Then $G$ is oddly finite dimensional by
Theorem \ref{main1} + Theorem \ref{MorelTheorem}.
\end{pf}

\begin{proposition}
\label{b} Assume that $U$ is not projective, i.e. $V\neq \emptyset
$. Then $M(U)$ is finite dimensional and, for any decomposition
  $$
  M(U)=\mathbb Q\oplus H\; ,
  $$
the motive $H$ is oddly finite dimensional.
\end{proposition}

\begin{pf}
Since $V\neq \emptyset $, we have the canonical distinguished
triangle
  \begin{equation}
  \label{tochki1}
  M(V)\stackrel{}{\lra }M(W)\lra M^c(U)\lra M(V)[1]
  \end{equation}
in $\DM _{\mathbb Q}$, where the motive $M(V)$ is just a direct
sum ${\mathbb Q^{\oplus }}^n$ of $n$ copies of the unit motive
$\mathbb Q$. Let
  $$
  M(W)=\mathbb Q\oplus M^1(W)\oplus \mathbb Q(1)[2]
  $$
be a decomposition of the motive of the smooth projective curve
$W$ by means of a point from $V$, say $v\in V$. Then we rewrite
the triangle (\ref{tochki1}) as follows:
  \begin{equation}
  \label{tochki2}
  {\mathbb Q^{\oplus }}^n \stackrel{}{\lra }
  \mathbb Q\oplus M^1(W)\oplus \mathbb Q(1)[2]
  \lra M^c(U)\lra {\mathbb Q[1]^{\oplus }}^n
  \end{equation}

The composition
  $$
  \mathbb Q\to {\mathbb Q^{\oplus }}^n \lra
  \mathbb Q\oplus M^1(W)\oplus \mathbb Q(1)[2]
  \to \mathbb Q\; ,
  $$
corresponding to the point $v$, is an isomorphism. Applying Lemma
\ref{usefullemma} we have that (\ref{tochki2}) is a direct sum of
the distinguished triangles
  \begin{equation}
  \label{tochki3}
  {\mathbb Q^{\oplus }}^{n-1} \stackrel{}{\lra }
  M^1(W)\oplus \mathbb Q(1)[2]
  \lra M^c(U)\lra {\mathbb Q[1]^{\oplus }}^{n-1}
  \end{equation}
and
  $$
  \mathbb Q\stackrel{}{\lra }
  \mathbb Q
  \lra 0\lra \mathbb Q[1]
  $$
Since any map from $\mathbb Q$ into $\mathbb Q(1)[2]$ is zero, the
triangle (\ref{tochki3}) is a direct sum of the triangles
  \begin{equation}
  \label{tochki4}
  {\mathbb Q^{\oplus }}^{n-1} \stackrel{}{\lra }
  M^1(W)\lra N\lra {\mathbb Q[1]^{\oplus }}^{n-1}
  \end{equation}
and
  $$
  0\stackrel{}{\lra }
  \mathbb Q(1)[2]
  \lra \mathbb Q(1)[2]\lra 0\; ,
  $$
so that
  $$
  M^c(U)\cong N\oplus \mathbb Q(1)[2]\; .
  $$

Note that $M^1(W)$ is oddly finite dimensional by Th.
\ref{findimcurves}. Applying Theorem \ref{main1} jointly with
Theorem \ref{MorelTheorem} to the triangle (\ref{tochki4}) we
claim that the motive $N$ is oddly finite dimensional.

Furthermore, since $U$ is a smooth scheme of pure dimension one,
  $$
  M(U)\cong N^*(1)[2]\oplus \mathbb Q
  $$
by \cite{Voev}, Th. 4.3.7 (3), where $N^*$ is a motive dual to
$N$. The dualization is a tensor endofunctor of $\DM $
\footnote{recall that $\DM $ is rigid}, whence $N^*$ is oddly
finite dimensional because $N$ is so. The motive $N^*(1)[2]$ is
oddly finite dimensional as a product of motives with different
parities, see Th. \ref{tensprod}.

Assume now that we have a splitting $M(U)=H\oplus \mathbb Q$.
Since $M(U)$ is finite dimensional as a sum of finite dimensional
motives $N^*(1)[2]$ and $\mathbb Q$, one can use Proposition
\ref{uniqdecomp} to show that $H\cong N^*(1)[2]$.
\end{pf}

\begin{proposition}
\label{c} Assume that $X$ is not projective. Then $M(X)$ is finite
dimensional, and, moreover,
  $$
  M(X)=\mathbb Q\oplus D\; ,
  $$
where $D$ is oddly finite dimensional in $\DM _{\mathbb Q}$.
\end{proposition}

\begin{pf}
If $X$ is smooth then $U=X$ and the proposition follows from
Proposition \ref{b}. Assume that $X$ is singular. Again, for
simplicity, we will prove the proposition in the case when $X$ has
only one singular point. Due to Proposition \ref{b}, the present
proof is almost the same like the proof of Proposition \ref{a}. We
consider a blow up distinguished triangle
  \begin{equation}
  \label{2t}
  {\mathbb Q^{\oplus }}^n\lra \mathbb Q\oplus M(U)\lra M(X)
  \lra {\mathbb Q^{\oplus }}^n[1]
  \end{equation}
associated with the contraction $p$ of points $\{ u_1,\dots ,u_n\}
$ onto a singular point of $X$. By Lemma \ref{usefullemma} the
triangle (\ref{2t}) splits into two distinguished triangles
  \begin{equation}
  \label{3t}
  {\mathbb Q^{\oplus }}^{n-1}\stackrel{t}{\lra }M(U)\lra
  M(X)\lra {\mathbb Q[1]^{\oplus }}^{n-1}
  \end{equation}
and
  $$
  \mathbb Q\lra \mathbb Q\lra 0\lra \mathbb Q[1]\; ,
  $$
where the last one corresponds to, for example, the point $u_1$.
Let
  $$
  M(U)=\mathbb Q\oplus \tilde M(U)
  $$
be a splitting induced by a $k$-rational point on $U$, say, by any
point $u_i$. Then, as above, for any $x_i$, $i>1$, the composition
  $$
  \mathbb Q\to {\mathbb Q^{\oplus }}^{n-1}
  \stackrel{t}{\lra }M(U)=\mathbb Q\oplus \tilde M(U)\lra \mathbb Q
  $$
is equal to zero, whence the triangle (\ref{3t}) is a direct sum
of two distinguished triangles
  \begin{equation}
  \label{4t}
  {\mathbb Q^{\oplus }}^{n-1}\stackrel{t}{\lra }\tilde M(U)\lra
  D\lra {\mathbb Q[1]^{\oplus }}^{n-1}
  \end{equation}
and
  $$
  0\lra \mathbb Q\lra \mathbb Q\lra 0\; .
  $$
In particular,
  $$
  M(X)=\mathbb Q\oplus D\; .
  $$
Note that $\tilde M(U)$ is oddly finite dimensional by Prop.
\ref{b}. Then, again, applying Theorem \ref{main1} together with
Theorem \ref{MorelTheorem} to the triangle (\ref{4t}) we have that
$D$ is oddly finite dimensional because $\tilde M(U)$ is so.
\end{pf}

\subsection{\it A remark on surfaces}

Let $X$ be a smooth projective surface with
$p_g=q=0$\footnote{where $p_g$ is a geometrical genus and $q$ is
an irregularity of $X$} and let $b_2$ be its second Betti number
in the sense of some Weil cohomology theory. Assume that the
motive $M(X)$ is finite dimensional. Then it is evenly finite
dimensional, and, in fact, it is isomorphic to the direct sum of
$\mathbb Q$, $b_2$ copies of the Lefschetz motive $\mathbb
Q(1)[2]$ and its tensor square $\mathbb Q(2)[4]$, see \cite{GP1}.
Let $U$ be a Zariski open subset in $X$. Theorem \ref{main1} and
the methods from Section \ref{pfmain2} allow then to prove that
$M(U)$ is also finite dimensional. In particular, if Bloch's
conjecture holds for a smooth projective complex surface $X$ of
general type with $p_g=0$, i.e. $M(X)$ is finite dimensional, see
\cite{GP2}, then $M(U)$ is finite dimensional for any Zariski open
subset $U$ in $X$.

\begin{small}

\end{small}

\bigskip

\noindent {\tt guletskii@im.bas-net.by}

\noindent {\sc Institute of Mathematics, Surganova 11, 220072
Minsk, Belarus}




\end{document}